\documentclass{article}
\usepackage{amsthm,amsmath,amssymb}

\newcommand{\nc}{\newcommand}

\newtheorem{thm}{Theorem}[section]
\newtheorem{rmk}[thm]{Remark}

\newtheorem{prop}[thm]{Proposition}
\newtheorem{lemma}[thm]{Lemma}

\newtheorem{corollary}[thm]{Corollary}

\newtheorem{definition}[thm]{Definition}

\newenvironment{defin}{\begin{definition} \rm}{\end{definition}}
\newenvironment{cor}{\begin{corollary} \rm}{\end{corollary}}
\newenvironment{lem}{\begin{lemma}\rm }{\end{lemma}}

\nc{\Ext}{\operatorname{Ext}}
\nc{\NS}{\operatorname{NS}}
\nc{\Amp}{\operatorname{Amp}}
\nc{\Pic}{\operatorname{Pic}}
\nc{\Kom}{\operatorname{Kom}}
\nc{\Gr}{\operatorname{Gr}}
\nc{\Rep}{\operatorname{Rep}}
\nc{\Hom}{\operatorname{Hom}}
\nc{\RHom}{R\operatorname{Hom}}
\nc{\cRHom}{\operatorname{\mathcal{R}\mathcal{H}om}}
\nc{\cHom}{\operatorname{\mathcal{H}om}}
\nc{\End}{\operatorname{End}}
\nc{\Coh}{\operatorname{Coh}}
\nc{\Aut}{\operatorname{Aut}}
\nc{\Coker}{\operatorname{Coker}}
\nc{\coker}{\operatorname{coker}}
\nc{\Ker}{\operatorname{Ker}}
\nc{\img}{\operatorname{Im}}
\nc{\D}{\operatorname{D}}
\nc{\ch}{\operatorname{ch}}
\nc{\Stab}{\operatorname{Stab}}
\nc{\SL}{\operatorname{Stab}}
\nc{\rk}{\operatorname{rk}}
\nc{\GL}{\operatorname{GL}}
\nc{\Log}{\mathop{\mathrm{Log}}}
\nc{\abs}[1]{\lvert#1\rvert}
\nc{\Cone}{\operatorname{Cone}}
\nc{\id}{\operatorname{id}}

\nc{\cA}{{\mathcal A}}
\nc{\cB}{{\mathcal B}}
\nc{\cC}{{\mathcal C}}
\nc{\cD}{{\mathcal D}}
\nc{\cE}{{\mathcal E}}
\nc{\cF}{{\mathcal F}}
\nc{\cG}{{\mathcal G}}
\nc{\cH}{{\mathcal H}}
\nc{\cI}{{\mathcal I}}
\nc{\cJ}{{\mathcal J}}
\nc{\cK}{{\mathcal K}}
\nc{\cL}{{\mathcal L}}
\nc{\cM}{{\mathcal M}}
\nc{\cN}{{\mathcal N}}
\nc{\cO}{{\mathcal O}}
\nc{\cP}{{\mathcal P}}
\nc{\cQ}{{\mathcal Q}}
\nc{\cR}{{\mathcal R}}
\nc{\cS}{{\mathcal S}}
\nc{\cT}{{\mathcal T}}
\nc{\cU}{{\mathcal U}}
\nc{\cV}{{\mathcal V}}
\nc{\cW}{{\mathcal W}}
\nc{\cX}{{\mathcal X}}
\nc{\cY}{{\mathcal Y}}
\nc{\cZ}{{\mathcal Z}}

\nc{\bA}{{\mathbb A}}
\nc{\bB}{{\mathbb B}}
\nc{\bC}{{\mathbb C}}
\nc{\bD}{{\mathbb D}}
\nc{\bE}{{\mathbb E}}
\nc{\bF}{{\mathbb F}}
\nc{\bG}{{\mathbb G}}
\nc{\bH}{{\mathbb H}}
\nc{\bI}{{\mathbb I}}
\nc{\bJ}{{\mathbb J}}
\nc{\bK}{{\mathbb K}}
\nc{\bL}{{\mathbb L}}
\nc{\bM}{{\mathbb M}}
\nc{\bN}{{\mathbb N}}
\nc{\bO}{{\mathbb O}}
\nc{\bP}{{\mathbb P}}
\nc{\bQ}{{\mathbb Q}}
\nc{\bR}{{\mathbb R}}
\nc{\bS}{{\mathbb S}}
\nc{\bT}{{\mathbb T}}
\nc{\bU}{{\mathbb U}}
\nc{\bV}{{\mathbb V}}
\nc{\bW}{{\mathbb W}}
\nc{\bX}{{\mathbb X}}
\nc{\bY}{{\mathbb Y}}
\nc{\bZ}{{\mathbb Z}}

\begin{document}
\title{On stability manifolds of Calabi-Yau surfaces} \author{So
Okada\footnote{Address: Vivatsgasse 7, Bonn Germany 53111, Email:
okada@mpim-bonn.mpg.de}}

\maketitle
\begin{abstract}
 We prove some general statements on stability conditions of Calabi-Yau
 surfaces and discuss the stability manifold of the cotangent bundle of
 $\bP^{1}$.  Our primary interest is in {\it spherical} objects.
\end{abstract}
\section{Introduction}

The notion of stability conditions on triangulated categories was
formulated in \cite{BRD_ST}.  It organizes certain bounded
$t$-structures on a triangulated category into a complex manifold.  In
the case of Calabi-Yau spaces, this is expected to be an approximation
of the stringy K\"ahler moduli of $X$.

Stability conditions have been studied for one-dimensional spaces in
\cite{BRD_ST}, \cite{GKR}, \cite{Ok}, \cite{Ma1}, and \cite{BuKr},
higher-dimensional spaces in \cite{Th}, \cite{BRD_K3},
\cite{BRD_Local_CY}, \cite{BRD_KL}, \cite{BRD_FR}, \cite{Ma1},
\cite{Ma2}, \cite{AB}, \cite{To}, \cite{Hu}, \cite{Be}, and \cite{An},
and $A_{\infty}$-categories in \cite{Th}, \cite{Ta}, \cite{Wa}, and
\cite{KST}.  The stability manifold of the category $\cO$ for ${\it
sl}_{2}$ has been computed in \cite{Mi}.  Some general aspects have been
studied in \cite{AP} and \cite{GKR}.  The author recommends
\cite{BRD_ICM} and \cite{Do1}, \cite{Do2}, \cite{Do3} for an
introduction and the original physical motivation to this subject.

We begin with fundamental notions and properties of stability
conditions.  After preparation on spectral sequences and {\it
$n$-Calabi-Yau categories}, we will concentrate on stability conditions
on 2-Calabi-Yau categories. Our main result is the connectedness of the
stability manifold of the cotangent bundle of $\bP^{1}$.

\subsubsection*{Acknowledgment}
The author thanks his adviser I. Mirkovi\'c for helping him achieve his
goals in the graduate school and people who gaved him the opportunity to
spend a month at the 2005 Seattle AMS Summer Institute.  The author
thanks E. Markman, K.  Yoshioka, and T. Bridgeland for their
discussions.  This work was completed in May 06 (for a grant
application), except that it has been revised in response to the
referee's comments. The author thanks two groups A. Ishii, K. Ueda,
H. Uehara and D. Huybrechts, E. Macr\`i, P. Stellari for communicating
to him their related work in Aug 06.  The author thanks the
Max-Planck-Institut f\"{u}r Mathematik for their support and excellent
working conditions.

\subsection{Definitions}
Throughout this paper, $\cT_{0}$ is the bounded derived category of an
abelian category with enough injectives and $\cT$ is a full triangulated
subcategory of $\cT_{0}$. In addition, $\cT$ is assumed to be linear
over $\bC$ and of finite type; i.e., for objects $E,F\in \cT$,
$\Hom_{\cT}(E,F)$ is a vector space over $\bC$ and the vector space
$\oplus_{i}\Hom^{i}_{\cT}(E,F)$ is of finite dimension.

For example, $\cT$ can be the bounded derived category $D(X)$ of
coherent sheaves on a smooth projective variety $X$, and $\cT_{0}$ the
bounded derived category of quasi-coherent sheaves on $X$ (\cite[Section
II, Proposition 2.2.2]{BGI}).  Let $K(\cT)$ be the $K$-group of $\cT$.  For an
object $E\in \cT$, let $[E]$ be the class of $E$ in $K(\cT)$.

We will recall some notions from \cite{BRD_ST}.

\subsubsection{Stability conditions}
A {\it stability condition} $\sigma=(Z, \cP)$ on $\cT$ consists of a
group homomorphism $Z:K(\cT)\to \bC$, called the {\it central charge},
and a family $\cP(\phi)$, $\phi \in \bR$, of full abelian subcategories
of $\cT$, called the {\it slicing}.  These need to satisfy the following
conditions.  If for some $\phi\in \bR$, \ $E$ is a nonzero object in
$\cP(\phi)$, then for some $m(E)\in\bR_{>0}$, \
$Z(E)=m(E)\exp(i\pi\phi)$.  For each $\phi\in \bR$, \
$\cP(\phi+1)=\cP(\phi)[1]$. For real numbers $\phi_{1}>\phi_{2}$ and
objects $A_{i}\in \cP(\phi_{i})$, \ $\Hom_{\cT}(A_{1}, A_{2})=0$. For
any object $E \in \cT$, there exist real numbers $\phi_{1}> \cdots
>\phi_{n}$ and objects $H^{\phi_{i}}_{\sigma}(E)\in \cP(\phi_{i})$ such
that there exists a sequence of exact triangles $E_{i-1}\to E_{i}\to
H_{\sigma}^{\phi_{i}}(E)$ with $E_{0}=0$ and $E_{n}=E$.  The sequence is
called the {\it Harder-Narasimhan filtration} (or {\it HN-filtration}
for short) of $E$. The HN-filtration of any object is unique up to
isomorphisms.

\subsubsection{Stability manifolds}
For an interval $I\subset \bR$, \ $\cP(I)$ denotes the smallest full
subcategory of $\cT$ that contains $\cP(\phi)$ for $\phi\in I$, it is
closed under extension; i.e., if $E\to G\to F$ is an exact triangle in
$\cT$ and $E, F\in \cP(I)$, then $G\in \cP(I)$. If the length of $I$ is
less than one, then $\cP(I)$ is a {\it quasi-abelian} category (in
particular, it is an exact category), whose exact sequences are
triangles of $\cT$ with vertices in $\cP(I)$.

A stability condition $\sigma=(Z, \cP)$ on $\cT$ is called {\it
locally-finite}, if for any $\phi\in \bR$, there exists a real number
$\eta >0$ such that $\cP((\phi-\eta,\phi+\eta))$ is of finite
length. The set of all locally-finite stability conditions on $\cT$ is
called the {\it stability manifold} of $\cT$ and denoted by
$\Stab(\cT)$.  The stability manifold of $\cT$ has a natural topology
and each connected component is a manifold locally modeled on some
topological vector subspace of $\Hom_{\bZ}(K(\cT), \bC)$.

\subsubsection{Some actions on stability manifolds}
Any stability manifold has a natural action of the group
$\widetilde{\GL^{+}}(2, \bR)$, the universal cover of
orientation-preserving transformations of $\GL(2, \bR)$.  In particular,
the group contains the following $\bC$-action for {\it rotation} and
{\it rescaling} of stability conditions; for $(Z, \cP)\in \Stab(\cT)$
and $z=x+i y\in \bC$, \ $z*(Z, \cP)$ is defined as $z*Z=e^{z}Z$ and
$(z*\cP)(\phi)=\cP(\phi-y/\pi)$ (\cite[Definition 2.3]{Ok}).

\subsubsection{Hearts of stability conditions} 
 For each $j\in \bR$, $\cP((j-1,j])$ and $\cP([j-1,j))$ are hearts of
 bounded $t$-structures.  By a {\it heart} of $\cT$, we mean the heart
 of any bounded $t$-structure of $\cT$. In particular, $\cP((0,1])$ is
 said to be the heart associated to a stability condition $\sigma=(Z,
 \cP)\in \Stab(\cT)$.  We will call all $c*\cP((0,1])$, $c\in \bC$,
 ``hearts of $\sigma$''.

\subsubsection{Semistable objects and stable objects}
 For a nonzero object $E\in \cP((j-1, j])$, the {\it phase} of $E$ is
 defined to be $\phi(E)=(1/\pi)\arg Z(E)\in (j-1, j]$. We say a real
 number $k$ is a {\it trivial phase} of an object $E\in \cT$, if
 $H^{k}_{\sigma}(E)$ is zero.

 For any $\phi\in \bR$, nonzero objects in $\cP(\phi)$ are called {\it
 semistable} objects.  For each object $E\in \cT$ and $k\in \bR$,
 $H^{k}_{\sigma}(E)$ is called the {\it semistable factor} of $E$ of
 the phase $k$. For each $k\in \bR$, any object $E\in \cP(k)$ has a
 Jordan-H\"older filtration in $\cP(k)$. A nonzero object $E\in\cP(k)$
 is called {\it stable} if it has no nontrivial subobject in $\cP(k)$.

\subsubsection{Jordan-H\"older blocks}
 \begin{defin}\label{def:JH}
  For an object $E\in \cT$, $k\in \bR$, and $\sigma\in \Stab(\cT)$, we
  will choose (non-canonical) ``{\it Jordan-H\"older blocks}'' (or {\it
  JH-blocks} for short) of $E$ denoted by $J_{\sigma}^{k}(E)$.  Let
  $A_{0}=0$ and $B_{0}=H^{k}_{\sigma}(E)$.  For $i>0$, let $A_{i}$ be a
  maximal subobject of $B_{i-1}$ such that all stable factors of $A_{i}$
  are isomorphic and let $B_{i}=B_{i-1}/A_{i}$. By the local-finiteness
  of $\sigma$, \ $B_{n}=0$ for some large enough $n$. We let
  $J_{\sigma}^{k}(E)=\{A_{1},\cdots, A_{n}\}$.
 \end{defin}

\begin{cor}\label{cor:JH}
 With the notation in Definition \ref{def:JH}, $\Hom_{\cT}(A_{i},
 B_{i})=0$ for any $0 \leq i \leq n$ and $\sum_{1\leq i\leq n}
 [A_{i}]=[H_{\sigma}^{k}(E)]$.
\end{cor}

  \section{Spectral sequences and $n$-Calabi-Yau categories}
  \subsection{Spectral sequences}
  For complexes $E, F\in \cT_{0}$ and a morphism of complexes $f:E\to
  F$, let $C(f)$ be the cone of $f$.  Let us say that $f$ is injective
  and splitting if $f$ is injective and it splits in each degree.

  \begin{lem}\label{lem:fund}
   Let $n\in \bZ_{> 0}$.  For each $0 \leq i \leq n$, let $F_{i}\in
   \cT_{0}$ be a complex.  For each $0 \leq i < n$, let $f_{i}$ be a
   morphism of complexes $f_{i}:F_{i} \to F_{i+1}$. Then there exist
   complexes of injective objects $\tilde{F}_{i}\in \cT_{0}$ and
   morphisms of complexes $\tilde{f}_i:{F}_{i}\to \tilde{F}_{i+1}$ with
   the following properties: morphisms
   $F_{i}\stackrel{f_{i}}{\to}F_{i+1}$ and
   $\tilde{F}_{i}\stackrel{\tilde{f}_{i}}{\to}\tilde{F}_{i+1}$ are
   isomorphic: $\tilde{f}_{i}$ is injective and splitting.
   \end{lem}

  \begin{proof}
   For $F_{0}$ and $F_{1}$, choose quasi-isomorphic complexes of
   injective objects $\dot{F_{0}}$ and $\dot{F_{1}}$.  Then
   $f_{0}:F_{0}\to F_{1}$ is isomorphic to a morphism of complexes
   $\dot{f_{0}}:\dot{F_{0}}\to \dot{F_{1}}$.  Let
   $\tilde{F_{0}}=\dot{F_{0}}$, $\tilde{F_{1}}=\dot{F_{1}}\oplus
   C(\dot{f_{0}})$, and $\tilde{f_{0}}=\dot{f_{0}}\oplus \alpha_{0}$,
   where $\alpha_{0}:\dot{F_{1}}\to C(\dot{f_{0}})$ is the canonical
   morphism. Then $f_{0},\dot{f_{0}}$, and $\tilde{f_{0}}$ are
   isomorphic in $\cT_{0}$, and $\tilde{f_{0}}$ is injective in each
   degree, since $\alpha_{0}$ is injective in each degree.  Moreover,
   $\tilde{f_{0}}$ splits in each degree, since $\tilde{F_{0}}$ and $
   \tilde{F_{1}}$ are injective objects.  Now, proceed by induction.
  \end{proof}
  
  For a heart $\cA$ of $\cT$, let $\tau^{\cA}$ denote the truncation
  functor.  For any object $E\in \cT$, by Lemma \ref{lem:fund}, the
  sequence of canonical morphisms from $\tau_{\leq i-1}^{\cA}(E)$ to
  $\tau_{\leq i}^{\cA}(E)$ can be realized as a sequence of injective
  and splitting morphisms of complexes $\tilde{\tau}_{\leq i}^{\cA}(E)$.

\begin{defin}\label{def:spec}
 For an object $E\in \cT$ and integers $i\leq j$, let
 $\tilde{\tau}_{(i,j]}^{\cA}(E)=\tilde{\tau}^{\cA}_{\leq j}(E)/
 \tilde{\tau}^{\cA}_{\leq i}(E)$, in particular,
 $\tilde{\tau}^{\cA}_{(i-1,i]}[i](E)\cong H_{\cA}^{i}(E)$ in $\cT_{0}$.
 For an object $E\in\cT$ and each $i\in \bZ$, let $e^{\cA}_{i}(E)$ be
 the connecting morphism from $H^{i}_{\cA}(E)$ to $H^{i-1}_{\cA}(E)[2]$,
 in the exact triangle $\tilde{\tau}_{(i-2,i-1]}^{\cA}(E)
 \to\tilde{\tau}_{(i-2,i]}^{\cA}(E) \to
 \tilde{\tau}_{(i-1,i]}^{\cA}(E)$, here,
 \begin{align*}
  H^{i}_{\cA}(E)[-i] 
  \cong \tilde{\tau}_{(i-1,i]}^{\cA}(E)
  &\stackrel{e^{\cA}_{i}(E)[-i]}{\longrightarrow}
  H^{i-1}_{\cA}(E)[2-i]
  \cong  \tilde{\tau}_{(i-2,i-1]}^{\cA}(E)[1].
 \end{align*}
\end{defin}

\begin{defin}
 For objects $E,F\in \cT$, $p,q\in \bZ$, and a heart $\cA$ of $\cT$, let
 \begin{align*}
  E_{2,\cA}^{p,q} (E,F) &=
  \oplus_{i\in \bZ}\Hom_{\cT}^{p}(H^{i}_{\cA}(E), H^{i+q}_{\cA}(F)).
 \end{align*}
 For $\oplus f_{i}$ in $\oplus_{i\in \bZ}\Hom_{\cT}^{p}(H^{i}_{\cA}(E),
 H^{i+q}_{\cA}(F))=E_{2, \cA}^{p,q}(E,F)$, let
\begin{align*}
 d_{2, \cA}^{p,q}(E,F)(\oplus f_{i})
 &= \oplus_{i\in \bZ}((-1)^{p+q}f_{i-1}\circ
 e^{\cA}_{i}(E)-e^{\cA}_{i+q}(F)\circ f_{i})\in E_{2, \cA}^{p+2,q-1}(E,F).
\end{align*}
\end{defin}

\begin{prop}\label{prop:spec}
 For any heart $\cA$ of $\cT$ and any objects $E,F\in \cT$, there exists
 a spectral sequence converging to $\Hom_{\cT}^{n}(E,F)$ with its
 $(p,q)$-components and differentials on the second sheet given by
 $E_{2,\cA}^{p,q}(E,F)$ and $d_{2, \cA}^{p,q}(E,F)$.
\end{prop}
\begin{proof}
 For $P=E$ or $P=F$, we define a decreasing finite splitting sequence of
 subcomplexes $\cF^{i}( P)=\tilde{\tau}^{\cA}_{\leq -i} (P)$.  By
 \cite[3.1.3.4]{BBD}, applied to $\cT_{0}$, there exists a spectral
 sequence $E_{1}^{pq}=\oplus_{j-i=p}\Hom_{\cT}^{p+q}(\Gr^{i}(E),
 \Gr^{j}(F))$ that converges to $\Hom_{\cT}^{p+q}(E,F)$.  Here,
 $\Gr^{i}(E)=\cF^{i}(E)/\cF^{i+1}(E)\cong H_{\cA}^{-i}(E)[i]$.  With new
 variables $q'=-p$, $p'=2p+q$, $i'=-i$, and $j'=-j$, \ $E_{1}^{pq}$
 reads
 $\oplus_{q'+i'=j'}\Hom_{\cT}^{p'}(H_{\cA}^{i'}(E),H_{\cA}^{j'}(F))$.
 For $n\in \bZ_{>0}$, $(p,q)\mapsto (p+n,q-n+1)$ translates into
 $(p',q')\mapsto (p'+n+1, q'-n)$.  Observe that because of change of
 variables, term $E_{n}$ in the spectral sequence from
 \cite[3.1.3.4]{BBD} is now viewed as $E_{n+1}$.
\end{proof}

\begin{lem}\label{lem:id}
 For a heart $\cA$ of $\cT$ and an object $E\in \cT$, let $\id_{i}$ be
 the identity morphism of $H_{\cA}^{i}(E)$.  If $E$ is not a zero
 object, then $\Ker d_{2, \cA}^{0,0}(E,E)$ contains the one-dimensional
 vector space $(\oplus_{i}\id_{i})\otimes \bC$.
\end{lem}
\begin{proof}
 Here, $\oplus_{i}\id_{i}\in E_{2,\cA}^{0,0}(E,E)$ and
 $d_{2,\cA}^{0,0}(E,E)(\oplus_{i}\id_{i})=0$.
\end{proof}

\subsection{$n$-Calabi-Yau categories}

The dual of a vector space $V$ will be written $V^{*}$.
\begin{defin}\cite[Definition 3.1]{BoKa}
 A covariant additive functor $S:\cT\to \cT$ that commutes with shifts
 is called a {\it Serre functor}, if it is a category equivalence, and
 for any objects $E,F\in \cT$, there exist bi-functorial isomorphisms
 $\psi_{E,F}:\Hom_{\cT}(E,F)\cong \Hom_{\cT}(F,S(E))^{*}$ such that the
 composite $(\psi^{-1}_{S(E), S(F)})^{*}\circ
 \psi_{E,F}:\Hom_{\cT}(E,F)\to
 \Hom_{\cT}(F,S(E))^{*}\to\Hom_{\cT}(S(E),S(F))$ coincides with the
 isomorphism induced by $S$.
\end{defin}

By \cite[Proposition 3.4 b]{BoKa}, a Serre functor of $\cT$, if it exists,
is unique up to a canonical isomorphism of functors.  We will call the
bi-functorial isomorphisms $\{\phi_{E,F}\}_{E,F\in \cT}$, the {\it Serre
duality} of $\cT$.

\begin{defin}\cite{Ko}
 A triangulated category $\cT$ is called an {\it $n$-Calabi-Yau
 category}, if the shift $[n]$ is the Serre functor.
\end{defin}

   \begin{defin}
    We define the {\it $\cT$-dimension} of a heart $\cA$ of $\cT$ as the
    supremum of $n$ such that $\Hom_{\cT}^{n}(E,F)\neq 0$ for objects
    $E,F\in \cA$.
   \end{defin}

   \begin{prop}\label{prop:dim}
    For any $n$-Calabi-Yau category $\cT$, the $\cT$-dimension of any
    heart of $\cT$ is $n$.
   \end{prop}
   \begin{proof}
    For a non-zero object $E\in\cA$, \ $\Hom_{\cT}^{n}(E,E)\cong
    \Hom_{\cT}(E,E)^{*} \neq 0$.  For $m>n$ and any objects $E,F\in
    \cA$, \ $\Hom_{\cT}^{m}(E,F) =\Hom_{\cT}(E,F[m])\cong
    \Hom_{\cT}(F[m],E[n])^{*}=\Hom_{\cT}^{n-m}(F,E)^{*}=0$.
   \end{proof}

\begin{cor}\label{cor:quotient}
 For an $n$-Calabi-Yau category $\cT$, a heart $\cA$ of $\cT$, and
 objects $E,F\in \cT$, if $p<0$ or $n<p$, then $E_{2, \cA}^{p,
 q}(E,F)=0$.
\end{cor}
\begin{proof}
 Here, $E_{2, \cA}^{p,q}(E,F)=\oplus_{i}\Hom_{\cT}^{p}(H_{\cA}^{i}(E),
 H_{\cA}^{i+q}(F))$, which is zero when $p<0$, since $\cA$ is a heart of
 $\cT$. When $p>n$, Proposition \ref{prop:dim} applies.
\end{proof}

\begin{prop}\label{prop:indec}
 For an $n$-Calabi-Yau category $\cT$, an object $E\in \cT$,
 and $\sigma \in \Stab(\cT)$, let $k_{1}> \cdots> k_{m}$ be all
 non-trivial phases of $E$.  If for some $s\in \bZ$, \ $k_{s-1}-k_{s}>
 n-1$, then $E$ is decomposable.
\end{prop}
\begin{proof}
 Here, $E$ is an extension $E'\to E \to E''$ with $E'$ (resp. $E''$)
 being an extension of $H^{k_{i}}_{\sigma}(E)$ for $i< s$ (resp. $s \leq
 i$).  For $i < s \leq j$, \  $\Hom_{\cT}^{1}(H^{k_{j}}_{\sigma}(E),
 H^{k_{i}}_{\sigma}(E))=0$. Hence, $\Hom_{\cT}^{1}(E'',E')=0$.
\end{proof}
For objects $E,F\in \cT$ and $i\in \bZ$, let $(E,F)^{i}=\dim
 \Hom_{\cT}^{i}(E,F)$.  The {\it Euler form} on $K(\cT)$ is
 $\chi(E,F)=\sum_{i\in \bZ}(-1)^{i}(E,F)^{i}$.  The quotient
 $N(\cT)=K(\cT)/K(\cT)^{\perp}$ is called the {\it numerical
 Grothendieck group}.  For an $n$-Calabi-Yau category $\cT$, the Euler
 form is (anti)symmetric depending on the parity of $n$ and factors
 through $N(\cT)$.

\section{Stability conditions on $2$-Calabi-Yau categories}
 From now on, $\cT$ is assumed to be $2$-Calabi-Yau.

\begin{lem}\label{lem:quotient}
 Let $\cA$ be a heart of $\cT$.  Then for any objects $E,F\in \cT$ and
 any $q\in \bZ$, \ $\Hom_{\cT}^{1+q}(E,F)$ has a filtration such that
 $\Ker d_{2, \cA}^{0,q+1}(E,F)$, $\Coker d_{2, \cA}^{0,q}(E,F)$, and
 $E_{2, \cA}^{1,q}(E,E)$ appear as distinct subquotients.
\end{lem}
\begin{proof}
 By Corollary \ref{cor:quotient}, differentials on the third sheet are
 zero.  Hence, $\Ker d_{2, \cA}^{0,q+1}(E,F)$ and $\Coker d_{2,
 \cA}^{0,q}(E,F)$ are subquotients of $\Hom_{\cT}^{1+q}(E,F)$. Since
 again by Corollary \ref{cor:quotient}, morphisms $d_{2, \cA}^{1,
 q}(E,E)$ and $d_{2, \cA}^{-1, q+1}(E,E)$ are zero, $E_{2,
 \cA}^{1,q}(E,E)$ is a subquotient of $\Hom_{\cT}^{1+q}(E,E)$.
\end{proof}

\begin{lem}\label{lem:Mukai}\cite[Lemma 5.2]{BRD_K3}
 Suppose $\cA$ is a heart of $\cT$ and $0 \to A \to B\to C \to 0$ is a
 short exact sequence in $\cA$ with $(A,C)^{0}=0$. Then
 $(A,A)^{1}+(C,C)^{1} \leq (B,B)^{1}$.
\end{lem}
\begin{proof}
 The equation $\chi([B],[B])=\chi([A]+[C], [A]+[C])$ reads
 $(B,B)^{1}=(A,A)^{1}+(C,C)^{1}+2((B,B)^{0}+(A,C)^{1}-
 ((A,A)^{0}+(C,A)^{0}+(C,C)^{0}))$. The inequality $(B,B)^{0}+(A,C)^{1}-
 ((A,A)^{0}+(C,A)^{0}+(C,C)^{0})\geq 0$ follows from the exact sequence
 $0\to \Hom_{\cT}(C,A)\to \End_{\cT}(B) \to \End_{\cT}(A)\oplus
 \End_{\cT}(C) \to \Hom^{1}_{\cT}(C,A)$, which is obtained from the
 condition $(A,C)^{0}=0$ and endomorphisms of the exact triangle $A\to B
 \to C$.
\end{proof}

\begin{lem}\label{lem:inequalities}
 For a heart $\cA$ of $\sigma \in \Stab(\cT)$ and an object $E\in \cT$,
 \begin{align*}
  (E,E)^{1} &\geq  \sum_{i\in \bZ}(H_{\cA}^{i}(E), H_{\cA}^{i}(E))^{1}
  \geq \sum_{k\in \bR}(H^{k}_{\sigma}(E), H^{k}_{\sigma}(E))^{1} 
  \geq  \sum_{k\in \bR, S\in J_{\sigma}^{k}(E)}(S,S)^{1}.
 \end{align*}
\end{lem}
\begin{proof}
 By Lemma \ref{lem:quotient}, $E_{2, \cA}^{1,0}(E,E)$ is a subquotient
 of $\Hom_{\cT}^{1}(E,E)$. Hence, $(E,E)^{1}\geq \dim
 E_{2,\cA}^{1,0}(E,E)=\sum_{i\in \bZ} (H^{i}_{\cA}(E),
 H^{i}_{\cA}(E))^{1}$.
 
 For some $j\in \bR$, \ $\cA=\cP((j-1, j])$.  So for each $i\in \bZ$ and
 $k\in (i+j-1,i+j]$, the HN-filtration of $H^{i}_{\cA}(E)$ gives a short
 exact sequence $0 \to A \to H^{i}_{\cA}(E)\to C \to 0$ in $\cA$ such
 that $A$ and $B$ are extensions of $H^{k'}_{\sigma}(E)$ for $k'>k$ and
 $k'\leq k$; in particular, $\Hom_{\cT}(A,B)=0$.  Hence, by Lemma
 \ref{lem:Mukai}, the second inequality follows.

 By Lemma \ref{lem:Mukai} and Corollary \ref{cor:JH}, the last
 inequality follows.
\end{proof}

\begin{defin}
 If an object $E\in \cT$ satisfies $\sum_{i}(E,E)^{i}=2$, then $E$ is
 called {\it spherical} (\cite[Definition 1.1]{ST}).
\end{defin}

\begin{lem}\label{cor:spherical}
 Let $\cA$ be a heart of $\cT$ and $E\in \cT$ be a spherical object.  If
 for some spherical object $S\in \cA$, every $H^{i}_{\cA}(E)$ is a
 multiple of $S$, then $E$ is a shift of $S$.
\end{lem}
\begin{proof}
 By taking a shift of $E$, for some $n\in \bZ_{\geq 0}$, we may suppose
 $H^{i}_{\cA}(E)$ is nonzero only for $0 \leq i \leq n$. Then, $E_{2,
 \cA}^{0,n+1}= \oplus_{i\in \bZ}\Hom_{\cT}(H^{i}_{\cA}(E),
 H^{i+n+1}_{\cA}(E))=0$. So, $\Coker d_{2, \cA}^{0,n+1}= E_{2,
 \cA}^{2,n}(E,E) = \Hom_{\cT}^{2}(H^{0}_{\cA}(E), H_{\cA}^{n}(E)) $,
 which is by Lemma \ref{lem:quotient}, a subquotient of
 $\Hom_{\cT}^{2+n}(E,E)$.  Since $E$ is spherical, $n=0$.
\end{proof}

\begin{lem}\label{lem:spherical_2}
 Let $\sigma\in \Stab(\cT)$, $E\in \cT$ be a non-semistable spherical
 object, and $k_{1}>\cdots>k_{n}$ be all nontrivial phases of $E$. If
 every $H^{k_{i}}_{\sigma}(E)$ is a multiple of a stable spherical
 object $S_{i}$, then $k_{s-1}-k_{s}<1$ for some $s$.
\end{lem}
\begin{proof}
 Since $E$ is spherical, it is indecomposable.  So by Proposition
 \ref{prop:indec}, $k_{i-1}-k_{i}\leq 1$.  If every $k_{i-1}-k_{i}=1$,
 then since $E$ is not semistable, by Lemma \ref{cor:spherical}, there
 exists $i$ such that $S_{i-1}[-1]\not\cong S_{i}$.  Since $S_{i-1}[-1]$
 and $S_{i}$ are non-isomorphic stable objects of the same phases,
 $\Hom_{\cT}^{1}(S_{i-1}, S_{i})= \Hom_{\cT}(S_{i-1}[-1], S_{i})=0$.  So
 $\Hom_{\cT}^{1}(H^{k_{i-1}}_{\sigma}(E), H^{k_{i}}_{\sigma}(E))=0$.
 Then, since for any $p < i \leq q$, we have
 $\Hom_{\cT}^{1}(H_{\sigma}^{k_{p}}(E),H_{\sigma}^{k_{q}}(E))=0$, $E$
 would be decomposable.
\end{proof}

\begin{rmk}\label{rmk:evenness}
 If $A\in \cT$ is stable for some stability condition, then $(A,A)^{1}$
 is even; because, the skew-symmetric, non-degenerate pairing
 $\Hom_{\cT}^{1}(A,A) \times \Hom_{\cT}^{1}(A,A) \to
 \Hom_{\cT}^{2}(A,A)\cong \Hom_{\cT}(A,A)^{*}=\bC$ implies
 $\Hom_{\cT}^{1}(A,A)$ is a symplectic vector space.
\end{rmk}

The pairing above is a simple case of the one in \cite[Proposition
I.1.4]{RV}.

\begin{defin}
 If an object $E\in \cT$ satisfies $(E,E)^{1}=0$, then $E$ is called
 {\it rigid} (\cite[Definition 3.1]{Mu}).
\end{defin}

\begin{lem}\label{lem:spherical_stable}
 For $\sigma\in \Stab(\cT)$, an object $E\in\cT$, and a
 nontrivial phase $k\in \bR$ of $E$, any rigid JH-block of
 $H^{k}_{\sigma}(E)$ is a multiple of a stable spherical object.
\end{lem}
\begin{proof}
 Let $S$ be a rigid JH-block of $H^{k}_{\sigma}(E)$.  Then $[S]=n[A]$
 for some stable object $A$ and $n> 0$.  Since $S$ is semistable and
 rigid, $\chi(S, S)= 2(S,S)^{0}$.  Since $A$ is stable, $0<
 \chi(S,S)=n^{2}(2(A,A)^{0} -(A,A)^{1})=n^{2}(2-(A,A)^{1})$. So, by
 Remark \ref{rmk:evenness}, the Euler form forces $(A,A)^{1}=0$.
\end{proof}

\begin{lem}\label{lem:decomp}
 For $\sigma\in \Stab(\cT)$, a rigid object $E\in \cT$, and a nontrivial
 phase $k$ of $E$, (a) any JH-block of $H^{k}_{\sigma}(E)$ is a multiple
 of a stable spherical object; and (b) if $J_{\sigma}^{k}(E)$ has more
 than one object, then there exist non-isomorphic stable spherical
 factors for $H^{k}_{\sigma}(E)$.
\end{lem}
\begin{proof}
 By Lemmas \ref{lem:inequalities} and \ref{lem:spherical_stable},
 (a) holds.  For (b), not all stable factors of $H_{\sigma}^{k}(E)$ are
 isomorphic; otherwise, $J_{\sigma}^{k}(E)$ would have only one object.
\end{proof}

\begin{prop}\label{prop:cutting_by_stability}
 For any $\sigma \in \Stab(\cT)$, if there exists a non-semistable
 spherical object, then in some heart of $\sigma$, there exist two
 non-isomorphic stable spherical objects.
\end{prop}
\begin{proof}
 Let $E\in \cT$ be a non-semistable spherical object and $k_{1}>\cdots
 >k_{n}$ be all nontrivial phases of $E$. Since $E$ is indecomposable,
 by Proposition \ref{prop:indec}, $k_{i-1}-k_{i} \leq 1$.

 If some $J_{\sigma}^{k_{i}}(E)$ has more than one object, then the
 statement follows by Lemma \ref{lem:decomp} (b).  Let us assume
 otherwise; by Lemma \ref{lem:decomp} (a), every $H^{k_{i}}_{\sigma}(E)$
 is a multiple of a stable spherical object.  Since $E$ is not
 semistable, Lemma \ref{lem:spherical_2} applies.
\end{proof}

\subsection{Twist functors}

   \begin{defin}
    For a spherical object $E\in \cT$ and an object $F\in\cT$, the cone of
    the evaluation map $\RHom_{\cT}(E, F)\otimes E \to F$ is denoted by
    $T_{E}(F)$, the {\it twist functor} of $E$ (\cite[Section 1.1]{ST}).
   \end{defin}
   
   For $\sigma=(Z, \cP)\in \Stab(\cT)$, a spherical object $E\in \cT$,
   and $T_{E}\in \Aut(\cT)$, let $T_{E}\sigma=T_{E}(\sigma)=
   ({T}_{E}{Z}, {T}_{E}{\cP}) =(Z\circ T_{E}^{-1}, T_{E}\circ \cP)$.

    \begin{lem}\label{lem:twist_nonsemi}
     For $\sigma=(Z, \cP)\in \Stab(\cT)$, let $E\in \cP(0)$ be a stable
     spherical object and $F\in \cP([0,1))$ be an object such that
     $\Hom_{\cT}(F,E)=0$. Then $F\in ({T}_{E}\cP)([0,1])$.
     \end{lem}
     \begin{proof}
      By the assumption on $F$, $\Hom^{i}_{\cT}(E,F)=0$ unless $i=0$ or
      $1$; so, $\RHom_{\cT}(E,F)\otimes E$ is an extension of multiples
      of $E$ and $E[-1]$.  Since $E=T_{E}(E[1])$ lies in
      $({T}_{E}\cP)(1)$, \ $\RHom_{\cT}(E,F)\otimes E$ lies in
      $({T}_{E}\cP)([0,1])$.  By the definition of ${T}_{E}\sigma$, \
      $T_{E}(F)\in ({T}_{E}\cP)([0,1))$.  Hence, the statement follows
      by the exact triangle $\RHom_{\cT}(E,F)\otimes E\to F\to
      T_{E}(F)$.
     \end{proof}
     \begin{cor}\label{cor:twist_nonsemi}
      In addition to the assumptions in Lemma \ref{lem:twist_nonsemi},
      assume further that $\phi(F)\in (0,1)$ and $\Hom_{\cT}(E,F)\neq
      0$.  Then $F$ is not semistable in ${T}_{E}{\sigma}$.
      \end{cor}
      \begin{proof}
       In this setting, $Z(F)$ is in the open upper-half plane of $\bC$.
       The same is then true for $(T_{E}Z)(F)=Z(F)-\chi(F,E)Z(E)$, since
       $Z(E)\in \bR$.  So if $F$ were semistable in ${T}_{E}{\sigma}$,
       then $F\in ({T}_{E}{\cP})([0,1])$ would imply that $F\in
       ({T}_{E}{\cP})((0,1))$.  However, $\Hom_{\cT}(E, F)\neq 0$ would
       contradict $({T}_{E}{\phi})(F)<({T}_{E}{\phi})(E)$.
    \end{proof}

    \begin{rmk}\label{rmk:assumptions}
     For $\sigma\in \Stab(\cT)$, we will refer to the following
     conditions: (a) there exist non-isomorphic stable spherical objects
     $E$ and $F$ in $\sigma$ such that $\Hom_{\cT}(E,F)\neq 0$ and
     $0=\phi(E)< \phi(F)<1$; (b) any two stable spherical objects of the
     same phases are isomorphic.
    \end{rmk}
    
    \begin{cor}\label{cor:twist_nonsemi_2}
     In addition to the assumptions in Corollary
     \ref{cor:twist_nonsemi}, assume that $F$ is spherical and the
     condition (b) in Remark \ref{rmk:assumptions}. Then,
     $H^{0}_{T_{E}\sigma}(F)=0$.
    \end{cor}
    \begin{proof}
     By the condition (b) in Remark \ref{rmk:assumptions},
     $E[-1]=T_{E}(E)$ is the only stable spherical object in
     $T_{E}\sigma$ of the phase zero.  So by Lemma \ref{lem:decomp} (a),
     $H^{0}_{T_{E}(\sigma)}(F)$ is a multiple of $E[-1]$.  Here, by
     Lemma \ref{lem:twist_nonsemi}, $F\in (T_{E}\cP)([0,1])$. So if
     $H^{0}_{T_{E}(\sigma)}(F)$ were not zero, then $\Hom_{\cT}(F,
     E[-1])$ would not be zero; however, since $E$ and $F$ are in some
     heart of $\sigma$, \ $\Hom_{\cT}(F, E[-1])=0$.
    \end{proof}
    
    \begin{prop}\label{prop:twist_nonsemi}
     For $\sigma=(Z, \cP)\in \Stab(\cT)$, assume the conditions (a) and
     (b) in Remark \ref{rmk:assumptions}.  Then for $E$ and $F$ in the
     condition (a) in Remark \ref{rmk:assumptions}, $F$ is not
     semistable in $T_{E}\sigma$ and lies in $(T_{E}\cP)((0,1])$.
    \end{prop}
    \begin{proof}
     Since $\phi(E)<\phi(F)$, \ $\Hom_{\cT}(F,E)=0$.  So by Lemma
     \ref{lem:twist_nonsemi}, $F\in (T_{E}\cP)([0,1])$.  Since
     $\phi(F)\in(0,1)$ and $\Hom_{\cT}(E,F)\neq 0$, by Corollary
     \ref{cor:twist_nonsemi}, $F$ is not semistable in $T_E \sigma$.
     Also, since $F$ is spherical and we are assuming the condition (b)
     in Remark \ref{rmk:assumptions}, by Corollary
     \ref{cor:twist_nonsemi_2}, $H^{0}_{T_{E}\sigma}(F)=0$.
    \end{proof}

  \section{Cotangent bundle of $\bP^{1}$}
  Let $Z=\bP^{1}$, $X$ the cotangent bundle of $Z$, $\Coh_{Z}(X)$ the
  category of the coherent sheaves of $X$ supported by $Z$, and $\cT$
  the full subcategory of $D(X)$ consisting of objects supported on $Z$.
  The space $X$ is the minimal resolution of the Kleinian singularity
  $\bC^{2}/\bZ_{2}$.   Let us prove the connectedness of $\Stab(\cT)$.

  \subsection{Pairs of stable spherical objects}

   \begin{lem}\label{lem:spherical_cot}
    For spherical objects $E,F\in\cT$, we have the following: (a) for
    some $s_{F}(E)=\pm 1$ and $p_{F}(E)\in \bZ$, \ $[E]=s_{F}(E)[F]+
    p_{F}(E)[\cO_{x}]$; (b) $p_{E}(F)=p_{E}(T_{E}(F))$.
   \end{lem}
   \begin{proof}
    We have $N(\cT)\cong \bZ \cdot Z(\cO_{Z})$, so
    $\chi(\cO_{Z},\cO_{Z})=2 =\chi(F,F)$ implies that $[F]$ is also a
    basis of $N(\cT)$.  Since $\bZ[\cO_{x}]=\Ker[K(\cT)\to N(\cT)]$, \
    $[F]$ and $[\cO_{x}]$ is a basis of $K(\cT)$.  Now, $s_{F}(E)=\pm
    1$, since $[E]=s_{F}(E) [F]$ in $N(\cT)$.  For the latter part,
    since $T_{E}[F]=[F]-\chi(E,F)[E]
    =(s_{E}-\chi(E,F))[E]+p_{E}(F)[\cO_{x}]$,
    $p_{E}(F)=p_{E}(T_{E}(F))$.
   \end{proof}

   \begin{defin}
    For a spherical object $E\in\cT$, we will call the sign
    $s_{\cO_{Z}}(E)$, the {\it sign of $E$}.
   \end{defin}   
   \begin{lem}\label{lem:compare_signs}
    If $E$, $F$, and $S$ are spherical objects in $\cT$ such that $E$
    and $F$ have different signs and $\phi(E[-1])<\phi(S)<
    \phi(F)<\phi(E)$, then $F$ and $S$ have the same signs.
   \end{lem}
    \begin{proof}
     Since $\phi(F)-\phi(S)\not\in \bZ$, \ $Z(F\oplus S)\neq 0$.  So if
     $F$ and $S$ had different signs, then by Lemma
     \ref{lem:spherical_cot} (a), $Z(F\oplus
     S)=(p_{E}(F)+p_{E}(S))Z(\cO_{x})\neq 0$.  Hence, $\phi(S)<\arg
     ((p_{E}(F)+p_{E}(S))Z(\cO_{x}))/\pi< \phi(F)$.  However, since
     $\phi(F)-\phi(E)\not\in\bZ$, again by Lemma \ref{lem:spherical_cot}
     (a), $Z(E\oplus F)= p_{E}(F) Z(\cO_{x}) \neq 0$, which implies
     $\phi(F)<\arg (p_{E}(F) Z(\cO_{x}))/\pi< \phi(E)$.
    \end{proof}

    By \cite[Theorem 1.2]{ST}, twist functors restrict to
    autoequivalences of $\cT$.

    \begin{lem}\label{lem:characterization_2}
     Suppose $\sigma\in \Stab(\cT)$ satisfies the conditions (a) and (b)
     in Remark \ref{rmk:assumptions}.  Then for $E$ and $F$ in the
     condition (a) in Remark \ref{rmk:assumptions}, if $E$ and $F$ have
     different signs, then there exist stable spherical objects $E'$ and
     $F'$ such that $E'$ and $F'$ have different signs,
     $\phi(E')<\phi(F')<\phi(E'[1])$, and $0<|(Z(E')+Z(F'))/Z(\cO_{x})|<
     |(Z(E)+Z(F))/Z(\cO_{x})|$.
    \end{lem}
    \begin{proof}
     By the conditions (a) and (b) in Remark \ref{rmk:assumptions}, and
     Proposition \ref{prop:twist_nonsemi}, $F$ is in $(T_{E}\cP)((0,1])$
     and not semistable in $T_{E}\sigma$.  Let $1=k_{1}>\cdots >k_{n}>0$
     be all nontrivial phases of $F$ in $T_{E}\sigma$.  By the condition
     (b) in Remark \ref{rmk:assumptions} and Lemma \ref{lem:decomp},
     each semistable factor of $F$ of the phase $k_{i}$ is a multiple of
     a stable spherical object $S_{i}$.

     Since $E[-1]$ and $F$ have the same signs and $(T_{E}\phi)(E[-1])
     <(T_{E}\phi)(S_{n})<(T_{E}\phi)(F)< (T_{E}\phi)(E)$, by Lemma
     \ref{lem:compare_signs}, $E[-1]$ and $S_{n}$ have the same
     signs. So $E'=T_{E}^{-1}(E)$ and $F'=T_{E}^{-1}(S_{n})$ have
     different signs.

     Since all $E[-1], S_{n}$, and $F$ have the same signs,
     $0=(T_{E}\phi)(E[-1]) <(T_{E}\phi)(S_{n})<(T_{E}\phi)(F)$ reads
     $0=\arg(Z[E])<\arg(Z([E]+p_{E}(F')[\cO_{x}])) <
     \arg(Z([E]+p_{E}(T_{E}^{-1}(F))[\cO_{x}]))$. Hence, $0
     <|p_{E}(F')|<|p_{E}(T_{E}^{-1}(F))|$.  By Lemma
     \ref{lem:spherical_cot} (b), $p_{E}(T_{E}^{-1}(F)) =p_{E}(F)$. So
     $0 <|p_{E}(F')|<|p_{E}(F)|$.

     Since $p_{E}(E')=p_{E}(E)=0$, $(Z(E')+Z(F'))/Z(\cO_{x})=p_{E}(F')$
     and $(Z(E)+Z(F))/Z(\cO_{x})=p_{E}(F)$. 
    \end{proof}

    \begin{prop}\label{prop:two_sphericals}
     For $\sigma \in \Stab(\cT)$, if there exists a non-semistable
     spherical object, then there exist two stable spherical objects in
     some heart of $\sigma$ such that they have no morphisms between
     them.
    \end{prop}
   \begin{proof}
    Since there exists a non-semistable spherical object, by Proposition
    \ref{prop:cutting_by_stability}, some heart of $\sigma$ contains two
    non-isomorphic stable spherical objects $E$ and $F$.

    Since any pair of non-isomorphic stable spherical objects of the
    same phases satisfies the conclusion, we may assume otherwise; i.e.,
    we may assume the condition (b) in Remark \ref{rmk:assumptions}.  In
    particular, $E$ and $F$ have different phases.

    By taking a shift of $E$ or $F$ if necessarily, we may assume $E$
    and $F$ have different signs. By using rotation and switching of $E$
    and $F$, we can assume $0=\phi(E)< \phi(F)$. Now, if $\Hom_{\cT}(E,
    F)=0$, then again $E$ and $F$ satisfy the conclusion.  So let
    $\Hom_{\cT}(E, F)\neq 0$, so for $E$ and $F$, the condition (a) in
    Remark \ref{rmk:assumptions} is satisfied.

    For our convenience, let $E_{0}=E$ and $F_{0}=F$.  By Lemma
    \ref{lem:characterization_2}, there exist stable spherical objects
    $E_{1}$ and $F_{1}$ such that $E_{1}$ and $F_{1}$ have different
    signs, $\phi(F_{1})<\phi(E_{1})<\phi(F_{1}[1])$, and
    $0<|(Z(E_{1})+Z(F_{1}))/Z(\cO_{x})|<|(Z(E_{0})+Z(F_{0}))/Z(\cO_{x})|$.
    
    So if we keep assuming $\Hom_{\cT}(E_{i},F_{i})\neq 0$, then we
    would obtain an infinite sequence of strictly decreasing positive
    integers $|p_{E_{i}}(F_{i})|$.  Hence, for some $E_{i}$ and $F_{i}$,
    \ $\Hom_{\cT}(E_{i}, F_{i})= \Hom_{\cT}(F_{i}, E_{i})=0$.
   \end{proof}

   \subsection{Pairs of stable spherical objects and autoequivalences}

   \begin{lem}\label{lem:computation}
    For $s,t,i\in \bZ$, $\Hom^{i}_{\cT}(\cO_{Z}(s),\cO_{Z}(t))= 0$ if
    and only if (a) $i=0$ and $s-t> 0$, (b) $i=1$ and $|s-t| < 2 $, or
    (c) $i=2$ and $s-t< 0$.
   \end{lem}
  \begin{proof}
   For (a), $\Hom_{\cT}(\cO_{Z}(s), \cO_{Z}(t)) =\Hom_{\cT}(\cO_{Z}
   (s-t), \cO_{Z})= 0$ if and only if $s-t> 0$; by the Serre duality,
   (c) follows.

   For (b), if $|s-t|<2$ and $k\neq 0$, then $\Hom_{D(Z)}^{k}(\cO,
   \cO(s-t))=\Hom_{D(Z)}^{k}(\cO(s-t), \cO)=0$.  So by \cite[Lemma
   4.6]{BRD_Local_CY}, $\Hom^{\bullet}_{\cT}(\cO_{Z}(s-t), \cO_{Z})
   =\Hom_{D(Z)}(\cO_{Z}(s-t), \cO_{Z}) \oplus \Hom_{D(Z)}(\cO_{Z},
   \cO_{Z}(s-t) )^{*}[-2]$; in particular, $\Hom^{1}_{\cT}(\cO_{Z}(s-t),
   \cO_{Z})=0$.
   
   For $|t-s| \geq 2$, let us recall that for the canonical bundle
   $\cO(-2)$ of $Z$, the spectral sequence in the proof of \cite[Lemma
   4.6]{BRD_Local_CY} reads $E_{2}^{p,q}=\Hom_{Z}^{p}(\cO,
   \cO(t-s)\otimes \wedge^{q} \cO(-2)) \Rightarrow
   \Hom_{\cT}^{p+q}(\cO_{Z}, \cO_{Z}(t-s))$.  Since the
   homological dimension of $\Coh(Z)$ is one, $E_{2}^{p,q}=0$ for $p>1$;
   also, since $\Coh(Z)$ is a heart of $D(\Coh Z)$, \ $E_{2}^{p,q}=0$ for
   $p<0$.

   If $t-s \geq 2$, then since $E_{2}^{0,1}=\Hom_{Z}(\cO,
   \cO(t-s-2))\neq 0$, \ $\Hom^{1}_{\cT}(\cO_{Z},\cO_{Z}(t-s))$ has a
   nonzero subquotient.  If $t-s \leq -2$, then  the Serre
   duality applies.
  \end{proof}

  For objects $E,F\in \cT$, we will write $H^{i}(E)$ and
  $E_{2}^{p,q}(E,F)$ for $H_{\Coh_{Z}(X)}^{i}(E)$ and
  $E_{2,\Coh_{Z}(X)}^{p,q}(E,F)$ (see Definition \ref{def:spec}).
  
  \begin{lem}\label{lem:difference}
   For objects $E,F\in \cT$ and $s,t,q\in \bZ$, let
   $\Hom_{\cT}^{i}(E,F)=0$ unless $i= 1$, $\cO_{Z}(s)$ is a summand of
   $H^{q}(E)$, and $F=\cO_{Z}(t)$.  Then, (a) $q\neq 0$ implies $|s-t|<
   2$, (b) $q=0$ and $\Hom_{\cT}(H^{-1}(E),F)=0$ implies $s-t<0$, and
   (c) $q=0$ and $\Hom_{\cT}^{2}(H^{1}(E), F)=0$ implies $s-t>0$.
  \end{lem}
  \begin{proof}
   For (a), if $|s-t|\geq 2$, then by Lemma \ref{lem:computation} (b),
   $\Hom_{\cT}^{1}(H^{q}(E), F)\neq 0$; i.e., $E_{2}^{1,-q}(E,F)\neq
   0$. Hence, by Lemma \ref{lem:quotient}, $\Hom_{\cT}^{1-q}(E,F)\neq
   0$. So $1-q=1$.

   For (b), by $E_{2}^{0,1}(E,F)=\Hom_{\cT}(H^{-1}(E),F)=0$, \ $\Coker
   d_{2}^{0,1}(E,F)$ is isomorphic to $E_{2}^{2,0}(E,F)$.  Since by
   Lemma \ref{lem:quotient}, $\Coker d_{2}^{0,1}(E,F)$ is a subquotient
   of $\Hom_{\cT}^{2}(E,F)=0$,
   $E_{2}^{2,0}(E,F)=\Hom_{\cT}^{2}(H^{0}(E), F)=0$.  So Lemma
   \ref{lem:computation} (c) applies.

   For (c), by $E_{2}^{2,-1}(E,F)=\Hom_{\cT}^{2}(H^{1}(E), F)=0$, \
   $\Ker d_{2}^{0,0}(E,F)$ is isomorphic to $E_{2}^{0,0}(E,F)$.  Since
   by Lemma \ref{lem:quotient}, $\Ker d_{2}^{0,0}(E,F)$ is a subquotient
   of $\Hom_{\cT}(E,F)=0$,
   $E_{2}^{0,0}(E,F)=\Hom_{\cT}(H^{0}(E),F)=0$. So Lemma
   \ref{lem:computation} (a) applies.
  \end{proof}

  \begin{rmk}\label{rmk:IU}
   By \cite[Section 5]{IU}, for any spherical object $E\in \cT$ and any
   $q\in \bZ$, we have some $v\in \bZ$ and $f_{q}, g_{q}\in \bZ_{\geq 0}$
   such that each $H^{q}(E)$ is isomorphic to $\cO_{Z}(v)^{f_{q}}\oplus
   \cO_{Z}(v-1)^{g_{q}}$, here, $l(E)$ is defined as
   $\sum_{q}(f_{q}+g_{q})$.
  \end{rmk}
  Let $S_{Z}(X)$ be the subgroup of $\Aut(\cT)$ generated by twists and
  shift functors on $\cT$.

   \begin{lem}\label{lem:tt}
    For spherical objects $E, F\in \cT$ and some $t\in \bZ$, if
    $l(E)>1$, $F$ is a shift of $\cO_{Z}(t)$, and for every irreducible
    summand $\cO_{Z}(s)$ of $H^{\bullet}(E)$ we have $|s-t|<2$, then for
    some $\Psi\in S_{Z}(X)$ we have $l(\Psi(E))<l(E)$ and
    $l(\Psi(F))=l(F)$.
   \end{lem}
  \begin{proof}
   By Remark \ref{rmk:IU}, for some $f_{q}, g_{q}\in \bZ_{\geq 0}$,
   either every $H^{q}(E)\cong \cO_{Z}(t+1)^{f_{q}}\oplus
   \cO_{Z}(t)^{g_{q}}$, or every $H^{q}(E)\cong \cO_{Z}(t)^{f_{q}}\oplus
   \cO_{Z}(t-1)^{g_{q}}$.  By \cite[Claim 5.2]{IU}, for the former case,
   $l(T_{\cO_{Z}(t)}(E))<l(E)$ and for the latter case,
   $l(T_{\cO_{Z}(t-1)}(E))<l(E)$.  Also for $F$, by \cite[Lemma 4.15
   (i)(1)]{IU}, \ $T_{\cO_{Z}(t)}(\cO_{Z}(t))=\cO_{Z}(t)[-1]$ and
   $T_{\cO_{Z}(t-1)}(\cO_{Z}(t))=\cO_{Z}(t-2)[1]$.
  \end{proof}

  \begin{lem}~\label{lem:two_sphericals_two} 
   For any $v\in \bZ$, we have $\Psi\in S_{Z}(X)$ such that
   $\{\Psi(\cO_{Z}(v)), \Psi(\cO_{Z}(v-1)[1])\}=\{\cO_{Z},
   \cO_{Z}(-1)[1]\}$.
  \end{lem}
  \begin{proof}
   By \cite[Lemma 4.15 (i)(2)]{IU}, $T_{\cO_{Z}(v)}\circ
   T_{\cO_{Z}(v-1)}$ acts as tensoring with $\cO_{Z}(-2)$.  So up to
   $S_{Z}(X)$, $v$ is zero or one.  
   For the latter case, by \cite[Lemma 4.15 (i)(1)]{IU},  
   $T_{\cO_{Z}}(\cO_{Z}[1])=\cO_{Z}$
   and  $T_{\cO_{Z}}(\cO_{Z}(1))=\cO_{Z}(-1)[1]$.
  \end{proof}

  \begin{prop}\label{prop:FM}
   For spherical objects $E,F\in \cT$ with $\Hom_{\cT}^{i}(E,F)=0$
   unless $i=1$, we have $\Psi\in S_{Z}(X)$ such that $\{\Psi(E),
   \Psi(F)\} =\{\cO_{Z}, \cO_{Z}(-1)[1]\}$.
  \end{prop}
  \begin{proof}
   By \cite[Proposition 5.1]{IU}, for some $\Psi_{1}\in S_{Z}(X)$ and
   $t\in \bZ$, \ $\Psi_{1}(F)=\cO_{Z}(t)$.  We will prove that for every
   summand $\cO_{Z}(s)$ of $H^{\bullet}(\Psi_{1}(E))$, \ $|s-t| < 2$.
   
   For some $q\neq 0$, if $\cO_{Z}(s)$ is a summand of
   $H^{q}(\Psi_{1}(E))$, then by Lemma \ref{lem:difference} (a),
   $|s-t|<2$.  Consider the case $q=0$.  If for a summand $\cO_{Z}(s)$
   of $H^{0}(\Phi_{1}(E))$, \ $s-t\geq 2$, then by Remark \ref{rmk:IU},
   any irreducible summand of $H^{-1}(\Psi_{1}(E))$ is $\cO_{Z}(r)$ of
   some $r-t>0$.  So by Lemma \ref{lem:computation} (a),
   $\Hom_{\cT}(H^{-1}(\Psi_{1}(E)),\Psi_{1}(F))=0$.  Here, by Lemma
   \ref{lem:difference} (b), $s-t<0$ in contradiction to $s-t\geq 2$.
   If for a summand $\cO_{Z}(s)$ of $H^{0}(\Phi_{1}(E))$, \ $s-t\leq
   -2$, then by Remark \ref{rmk:IU}, any irreducible summand of
   $H^{1}(\Psi_{1}(E))$ is $\cO_{Z}(r)$ of some $r-t<0$.  So by Lemma
   \ref{lem:computation} (c),
   $\Hom_{\cT}^{2}(H^{1}(\Psi_{1}(E)),\Psi_{1}(F))=0$.  Here, by Lemma
   \ref{lem:difference} (c), $s-t>0$ in contradiction to $s-t\leq -2$.
   
   Hence by Lemma \ref{lem:tt}, for some $\Psi_{2}$ and $l,m,n\in \bZ$,
   \ $\Psi_{2}\circ \Psi_{1}(E)=\cO_{Z}(m)[l]$ and $\Psi_{2}\circ
   \Psi_{1}(F)=\cO_{Z}(n)$.  By using the argument above on
   $\Psi_{2}\circ \Psi_{1}(E)$ and $\Psi_{2}\circ \Psi_{1}(F)$, \
   $|m-n|<2$. Furthermore, we will prove that either $m=n-1$ and $l=1$,
   or $n=m-1$ and $l=-1$.
   
   If $m=n-1$, then by Lemma \ref{lem:computation} (a),
   $\Hom_{\cT}^{l}(\Phi_{2}\circ \Phi_{1}(E),\Phi_{2}\circ
   \Phi_{1}(F))\neq 0$, which by the assumption, implies $l=1$.  If
   $n=m-1$, then by Lemma \ref{lem:computation} (c), $\Hom_{\cT}^{l+2}
   (\Phi_{2}\circ \Phi_{1}(E),\Phi_{2}\circ \Phi_{1}(F))\neq 0$, which
   by the assumption, implies $l+2=1$. If $m$ were equal to $n$, then by
   Lemmas \ref{lem:computation} (a) and \ref{lem:computation} (c), both
   $\Hom_{\cT}^{l}(\Phi_{2}\circ \Phi_{1}(E),\Phi_{2}\circ \Phi_{1}(F))$
   and $\Hom_{\cT}^{l+2}(\Phi_{2}\circ \Phi_{1}(E),\Phi_{2}\circ
   \Phi_{1}(F))$ would be nonzero, which contradicts the assumption.

   So, when $l=1$, $\Phi_{2}\circ \Phi_{1}(E)=\cO_{Z}(n-1)[1]$ and
   $\Phi_{2}\circ \Phi_{1}(F)=\cO_{Z}(n)$; also, when $l=-1$,
   $\Phi_{2}\circ \Phi_{1}(E)[1]=\cO_{Z}(m)$ and $\Phi_{2}\circ
   \Phi_{1}(F)[1]=\cO_{Z}(m-1)[1]$.  Now, Lemma
   \ref{lem:two_sphericals_two} applies.
  \end{proof}

\subsection{Connectedness}
 By \cite[Theorem 1.3]{BRD_KL}, a connected component
 $\Stab_{0}(\cT)\subset \Stab(\cT)$ is invariant under $S_{Z}(X)$.  Let
 $\cA\subset \cT$ be the smallest extension-closed full subcategory
 containing $\cO_{Z}$ and $\cO_{Z}(-1)[1]$, and $U$ be a subset of
 $\Stab(\cT)$ consisting of stability conditions $(Z,\cP)$ such that
 $\cP((0,1])=\cA$, and $\img Z(\cO_{Z})$, $\img Z(\cO_{Z}(-1)[1])>0$
 (\cite[Lemma 3.1]{BRD_KL}).  By \cite[Lemma 3.6]{BRD_KL}, for any
 $\sigma \in \Stab_{0}(\cT)$, there exists $\Psi \in S_{Z}(X)$ such that
 some rotation of $\Psi(\sigma)$ lies in the closure of $U$.
 
 Hence, to show the connectedness of $\Stab(\cT)$, we will prove that
 for any $\sigma \in \Stab(\cT)$, there exists $\Psi\in S_{Z}(X)$ such
 that $\cA$ is a heart of $\Psi(\sigma)$.

\begin{thm}
 $\Stab(\cT)$ is connected.
\end{thm}
\begin{proof}
 Let $\sigma\in \Stab(\cT)$.  Suppose all spherical objects are
 semistable.  Then, by Lemma \ref{lem:computation} (a), for any $v\in
 \bZ$, \ $\phi(\cO_{Z}(v-1))\leq \phi(\cO_{Z}(v)) $.  We will prove that
 for some $w\in \bZ$, \ $\phi(\cO_{Z}(w-1))< \phi(\cO_{Z}(w))$.

 Let us assume otherwise; i.e., for any $v\in \bZ$,
 $\phi(\cO_{Z}(v-1))=\phi(\cO_{Z}(v))$.  Then, since by Lemma
 \ref{lem:computation}, $\RHom_{\cT}(\cO_{Z}(1), \cO_{Z})\otimes
 \cO_{Z}(1)$ is a multiple of $\cO_{Z}(1)[-2]$, the exact triangle
 $\cO_{Z}\to T_{\cO_{Z}(1)}(\cO_{Z}) \to \RHom_{\cT}(\cO_{Z}(1),
 \cO_{Z})\otimes \cO_{Z}(1)[1]$ would be the nontrivial HN-filtration of
 the semistable spherical object $T_{\cO_{Z}(1)}(\cO_{Z})$.

 So, since by Lemma \ref{lem:computation} (b),
 $\phi(\cO_{Z}(w+1)[-1])\leq \phi(\cO_{Z}(w-1))$, \
 $\phi(\cO_{Z}(w+1)[-1])\leq \phi(\cO_{Z}(w-1))<\phi(\cO_{Z}(w)) \leq
 \phi(\cO_{Z}(w+1))$. Hence, $O_{Z}(w)$ and $\cO_{Z}(w-1)[1]$ are in a
 heart of $\sigma$.  Now, Lemma \ref{lem:two_sphericals_two} applies.

 If not all spherical objects are semistable, then by Proposition
 \ref{prop:two_sphericals}, there exist non-isomorphic stable spherical
 objects $E$ and $F$ in some heart of $\sigma$ such that
 $\Hom_{\cT}(E,F)=\Hom_{\cT}(F,E)=0$.  By the Serre duality,
 $\Hom_{\cT}^{2}(E,F)=0$. Since $E$ and $F$ are in some heart of
 $\sigma$, for $i<0$, \ $\Hom_{\cT}^{i}(E,F)=0$.  So by Proposition
 \ref{prop:dim}, $\Hom_{\cT}^{i}(E,F)=0$ unless $i=1$.  Now, Proposition
 \ref{prop:FM} applies.
\end{proof}

\end{document}